\documentclass[11pt]{amsart}

\usepackage{amsfonts,amssymb,stmaryrd,amscd,amsmath,latexsym,amsbsy}

\usepackage{amssymb}
\usepackage{amsfonts}
\usepackage{latexsym}

\newtheorem{theorem}{Theorem}[section]

\newtheorem{proposition}[theorem]{Proposition}

\theoremstyle{definition}

\newtheorem{conjecture}[theorem]{Conjecture}

\newcommand{\ben}{\begin{enumerate}}
\newcommand{\een}{\end{enumerate}}

\theoremstyle{plain}

\newtheorem*{sol}{Solution}

\theoremstyle{definition}

\theoremstyle{remark}

\newcommand{\solu}[1]{\begin{sol}{\bf (\ref{#1})}}

\begin{document}

\title{Casimirs of the Goldman Lie algebra of a closed surface}

\author{Pavel Etingof}
\address{Department of Mathematics, Massachusetts Institute of Technology,
Cambridge, MA 02139, USA}
\email{etingof@math.mit.edu}

\maketitle

\section{Introduction}

Let $\Sigma$ be a connected closed oriented surface of genus $g$.
In 1986 Goldman \cite{Go} attached to $\Sigma$ a Lie algebra 
$L=L(\Sigma)$, later shown by Turaev (\cite{Tu}) to have a natural structure 
of a Lie bialgebra. It is defined as follows. As a vector space,
$L$ has a basis $e_\gamma$ labeled by conjugacy classes $\gamma$ 
in the fundamental group $\pi_1(\Sigma)$, geometrically represented
by closed oriented curves on $\Sigma$ without a base point. 
To define the commutator $[e_{\gamma_1},e_{\gamma_2}]$, 
one needs to bring the two curves $\gamma_1,\gamma_2$ into 
general position by isotopy, and then for each intersection point 
$p_i$ of the two curves, define $\gamma_{3i}$ 
to be the curve obtained by tracing $\gamma_1$ and then $\gamma_2$ 
starting and ending at $p_i$. Then one defines 
$[e_{\gamma_1},e_{\gamma_2}]$ to be 
$\sum_i \varepsilon_i e_{\gamma_{3i}}$, where 
$\varepsilon_i=1$ if $\gamma_1$ approaches $\gamma_2$ from the 
right at $p_i$ (with respect to the orientation of $\Sigma$), 
and $-1$ otherwise. 

The combinatorial structure of $L$ has been much studied; see
e.g. \cite{C,Tu}. 
However, many problems about the structure of $L$ remained open. 
In particular, in 2001, M. Chas and D. Sullivan communicated 
to me the following conjecture. 

\begin{conjecture}
The center of $L$ is spanned by 
the element $e_1$, where $1\in \pi_1(\Sigma)$ is the trivial
loop. 
\end{conjecture}

In this paper, we will prove this conjecture. 
In fact, we prove a more general result. 

\begin{theorem}\label{cs}
The Poisson center of the Poisson algebra $S^\bullet L$ is $Z=\Bbb C[e_1]$.   
\end{theorem}

The proof of the theorem occupies the rest of the paper. 

{\bf Remark.} A quiver theoretical analog of 
Theorem \ref{cs} is given in \cite{CEG}. 
It claims that if $\Pi$ is the preprojective algebra
of a quiver $Q$ which is not Dynkin or affine Dynkin,
then the Poisson center of $S^\bullet L$ (where $L=\Pi/[\Pi,\Pi]$ 
is the necklace Lie algebra attached to $\Pi$) consists of polynomials 
in the vertex idempotents. 

\section{Proof of the theorem}

\subsection{Moduli spaces of flat bundles}
We will assume that $g>1$, since 
in the case $g\le 1$ the theorem is easy. 

Recall that the fundamental group $\Gamma=\pi_1(\Sigma)$ is generated by 
$X_1,...,X_g$, $Y_1,...,Y_g$ with defining relation 
\begin{equation}\label{rel}
\prod_{i=1}^g X_iY_iX_i^{-1}Y_i^{-1}=1.
\end{equation}

Thus we can define the scheme of homomorphisms 
$\widetilde M_g(N)={\rm Hom}(\Gamma,GL_N(\Bbb C))$
to be the closed subscheme in $GL_N(\Bbb C)^{2g}$ defined by equation 
(\ref{rel}). One can also define the moduli scheme of representations
(or equivalently, of flat connections on $\Sigma$) 
to be the categorical quotient 
$M_g(N)=\widetilde M_g(N)/PGL_N(\Bbb C)$.

The schemes $\widetilde M_g(N)$ and $M_g(N)$ carry the 
Atiyah-Bott Poisson structure; its algebraic presentation may be found 
in \cite{FR} (using r-matrices) and \cite{AMM} (using quasi-Hamiltonian 
reduction); see also \cite{Go}. 

Let us recall the following known results about these schemes,
which we will use in the sequel. 

\begin{theorem}\label{modsp}
(i) $\widetilde M_g(N)$ and $M_g(N)$ are reduced. 

(ii) $\widetilde M_g(N)$ is a complete intersection
in $GL_N(\Bbb C)^{2g}$.

(iii) $\widetilde M_g(N)$ and $M_g(N)$ are irreducible algebraic varieties.
Their generic points correspond to irreducible representations 
of $\Gamma$.

(iv) The Poisson structure on $M_g(N)$ is generically symplectic.
\end{theorem}

\begin{proof}
Let $\widetilde M'_g(N)$ be the algebraic variety 
corresponding to the scheme $\widetilde M_g(N)$. It is shown in \cite{Li} 
that this variety is irreducible. Moreover, it is clear that 
the generic point of this variety corresponds to an irreducible representation 
of $\Gamma$ (we can choose $X_i,Y_i$ generically for $i<g$ and then solve 
for $X_g,Y_g$). It is easy to show that 
near such a point the map $\mu: GL(N)^{2g}\to SL(N)$ 
given by the left hand side of (\ref{rel}) is a submersion. 
This implies (ii). We also see that $\widetilde M_g(N)$ is 
generically reduced. Since it is a complete intesection, it is Cohen-Macaulay 
and hence reduced everywhere. Thus we get (i) and (iii). 
Property (iv) is well known and is readily seen from \cite{FR} or \cite{AMM}. 
The theorem is proved. 
\end{proof} 

\subsection{Injectivity of the Goldman homomorphism}

Now let us return to the study of the Lie algebra $L$. 
To put ourselves in an algebraic framework, we note 
that $L$ is naturally identified with $A/[A,A]$, where $A=\Bbb C[\Gamma]$ 
is the group algebra of $\Gamma$. Thus, elements of $L$ can be represented by 
linear combinations of cyclic words in $X_i^{\pm 1},Y_i^{\pm 1}$. 

In \cite{Go}, Goldman defined a homomorphism of Poisson algebras
$$
\phi_N: S^\bullet L\to \Bbb C[M_g(N)]
$$
defined by the formula $\phi_N(w)(\rho)={\rm Tr}(\rho(w))$, 
where $\rho$ is an $N$-dimensional representation of $\Gamma$
and $w$ is any cyclic word representing an element of $L$. 
It follows from Weyl's fundamental theorem of invariant theory 
that the Goldman homomorphism is surjective. 

Let $L_+\subset L$ be the linear span of 
the elements $e_\gamma-e_1$. Obviously, 
we have $L=L_+\oplus \Bbb C e_1$, 

\begin{proposition} \label{inj} 
For any finite dimensional subspace $Y\subset S^\bullet L_+$, 
there exists an integer $N(Y)$ such that for $N\ge N(Y)$, 
the map $\phi_N|_Y$ is injective. 
\end{proposition}

\begin{proof}
Let $K(N)$ be the kernel of $\phi_N$ on $S^\bullet L_+$. 
It is clear that $K(N+1)\subset K(N)$ 
(as $\phi_{N+1}(e_\gamma-e_1)(\rho\oplus \Bbb C)=
\phi_{N}(e_\gamma-e_1)(\rho)$). Thus it suffices to show that 
$\cap_{N\ge 1}K(N)=0$.

Assume the contrary. Then there exists an element $0\ne f\in S^\bullet
L_+$ such that $\phi_N(f)=0$ for all $N$. 

Recall that according to \cite{FiR},
the group $\Gamma$ is {\bf conjugacy separable}, i.e., 
if elements $\gamma_0,...,\gamma_m$ are pairwise not
conjugate in $\Gamma$
then there exists a finite quotient $\Gamma'$ of $\Gamma$ such
that the images of $\gamma_0,...,\gamma_m$ are not conjugate in
$\Gamma'$. 

Now let $\gamma_0=1$ and 
$f=P(e_{\gamma_1}-e_1,...,e_{\gamma_m}-e_1)$, where $P$ is
some polynomial. 
Let $\Gamma'$ be the finite group as above, 
$V_1,...,V_s$ be the irreducible representations of $\Gamma'$,
and $\chi_1,...,\chi_s$ be their characters. Let $V=\oplus_j N_j
V_j$; we regard $V$ as a representation of $\Gamma$ and let $N=\dim V$. 
Then $\phi_N(f)(V)=P(w_1,...,w_m)$, where 
$w_i=\sum_j N_j(\chi_j(\gamma_i)-\chi_j(1))$. By representation
theory of finite groups, the matrix with entries 
$a_{ij}=\chi_j(\gamma_i)-\chi_j(1)$ has rank $m$; thus, 
there exist $N_j\ge 0$ such that $P(w_1,...,w_m)\ne 0$. 
For such $N_j$, $\phi_N(f)\ne 0$, which is a contradiction. 
\end{proof}

\subsection{Proof of Theorem \ref{cs}}

Now we are ready to prove Theorem \ref{cs}. 
Let $z$ be a central element of the Poisson algebra $S^\bullet L$. 
Consider the element $\phi_N(z)$. 
This is a regular function on $M_g(N)$
which Poisson commutes with all other functions (since $\phi_N$ 
is surjective). Since by Theorem \ref{modsp} 
the scheme $M_g(N)$ is in fact a variety, which is irreducible and 
generically symplectic, any Casimir on this variety must be a scalar. 

Since $S^\bullet L=S^\bullet L_+\otimes \Bbb C[e_1]$, we can write $z$ as 
$$
z=\zeta(e_1)+\sum_{j=1}^m \zeta_j(e_1)f_j, 
$$
were $f_j$ 
are linearly independent elements which belong
to the augmentation ideal of $S^\bullet L_+$, and 
$\zeta,\zeta_j\in \Bbb C[t]$. 
Applying $\phi_N$ to this equation, and using that
$\phi_N(e_1)=N$, 
we get that 
$$
\zeta(N)+\sum_{j=1}^m \zeta_j(N)\phi_N(f_j)=\gamma_N.
$$
Let $Y$ be the linear span of $1$ and $f_j$, $j=1,...,m$ in
$S^\bullet L_+$. By Proposition \ref{inj}, 
for $N\ge N(Y)$, we have 
$$
\zeta(N)+\sum_{j=1}^m \zeta_j(N)f_j=\gamma_N.
$$
Thus $\zeta_j(N)=0$ for $N\ge N(Y)$. Hence $\zeta_j=0$ for all
$j$ and $z=\zeta(e_1)$. The theorem is proved. 

{\bf Acknowledgements.} This work
was partially supported by the NSF grant DMS-9988796
and the CRDF grant RM1-2545-MO-03. 
I am grateful to M. Chas and D. Sullivan for posing the
problem, and to Hebrew University for hospitality. I also thank 
V. Ginzburg for useful discussions and M. Sapir for 
explanations and references about conjugacy separability.

\end{document}